\documentclass[10pt,reqno]{amsart}
\usepackage{amssymb}
\usepackage{amsfonts}
\usepackage[all]{xy}
\usepackage[english]{babel}

\begin{document}

\title{Distributional chaos and factors}

\author{Jana Dole\v zelov\'a-Hant\'akov\'a}

\address{J . Hant\'akov\'a, Mathematical Institute, Silesian University, CZ-746 01 
Opava, Czech Republic}

\email{jana.hantakova@math.slu.cz, }

\thanks{The research was supported by grant SGS/2/2013 from  the Silesian University in Opava. Support of this institution is gratefully acknowledged.}

\begin{abstract} 
 We show the existence of a dynamical system without any distributionally scrambled pair which is semiconjugated to a distributionally chaotic factor.\\
{\small {2000 {\it Mathematics Subject Classification.}}
Primary 37D45; 37B40.
}
\end{abstract}

\maketitle
\pagestyle{myheadings}
\markboth{Jana Dole\v zelov\'a-Hant\'akov\'a}
{Distributional chaos and factors}
\section{Introduction}
Semiconjugacy is used as a common tool for proving topological chaos or positive topological entropy. The usual
technique is to find a semiconjugacy $\pi$ with a chaotic system and transfer the chaos to the extension. By continuity of $\pi$, the
topological entropy of the extension is not smaller than the entropy
of factor system. Unfortunately, semiconjugacy may not automatically guarantee the distributional chaos, which was introduced in \cite{SchSm}. Authors in \cite{OW},\cite{FOW},\cite{W} developed several techniques for
proving distributional chaos via semiconjugacy, usually using a symbolic space as the factor space. Example in \cite{OW} shows the existence of distributionally chaotic factor which is semiconjugated to the system with no three-points distributionally scrambled sets. The aim of the paper is to improve this result and find a distributionally chaotic factor which has an extension without any distributionally scrambled pair.

\section{Terminology}
Let $(X,d)$ be a non-empty compact metric space.  Let us denote by $(X,f)$ the \emph{topological dynamical system}, where $f$ is a continuous self-map acting on $X$.  We define the \emph{forward orbit} of $x$, denoted by $Orb^+_f(x)$ as the set $\{f^n(x):n\geq0\}$.
Let $(X,f)$ and $(Y,g)$ be dynamical systems on compact metric spaces. A continuous map $\pi: X\rightarrow Y$ is called a semiconjugacy between f and g if $\pi$ is surjective and $\pi\circ f=g\circ\pi$. In this case we can say that $(Y,g)$ is a factor of the system $(X,f)$ or eqivalently $(X,f)$ is an extension of the system $(Y,g)$.

\paragraph{\bf Definition 1} A pair of two different points $(x_1, x_2)\in X^2$ is called \emph{scrambled} if 
\begin{equation}\liminf_{k\to\infty}d(f^k(x_i),f^{k}(x_j))=0\end{equation} and \begin{equation}\limsup_{k\to\infty} d(f^k(x_i),f^{k}(x_j))>0.\end{equation}
A subset $S$ of $X$ is called \emph{scrambled} if every pair of distinct points in $S$ scrambled. The system $(X,f)$ is called \emph{chaotic} if there exists an uncountable scrambled set.\\

\paragraph{\bf Definition 2} For a pair $(x_1, x_2)$ of
points in $X$, define the \emph{lower distribution function} generated by $f$ as
$$\Phi_{(x_1, x_2)}(\delta)=\displaystyle\liminf_{m\to\infty}\frac{1}{m}\#\{0<k<m;d(f^k(x_1),f^{k}(x_2))<\delta\},$$
and the \emph{upper distributional function} as 
$$\Phi^*_{(x_1, x_2)}(\delta)=\displaystyle\limsup_{m\to\infty}\frac{1}{m}\#\{0<k<m;d(f^k(x_1),f^{k}(x_2))<\delta\},$$
where $\#A$ denotes the cardinality of the set $A$.\\ 
A pair $(x_1, x_2)\in X^2$ is called 
\emph{distributionally scrambled of type 1} if 
$$\Phi^*_{(x_1, x_2)}\equiv 1 \mbox{  and  } \Phi_{(x_1, x_2)}(\delta)=0, \mbox{  for some  }0<\delta\le \text{diam }X ,
$$
\emph{distributionally scrambled of type 2} if $$\Phi^*_{(x_1, x_2)}\equiv 1 \mbox{  and  } \Phi_{(x_1, x_2)}<\Phi^*_{(x_1, x_2)},$$
\emph{distributionally scrambled of type 3} if $$\Phi_{(x_1, x_2)}<\Phi^*_{(x_1, x_2)}.
$$
The dynamical system $(X,f)$ is distributionally chaotic of type $i$ (DC$i$ for short), where $i=1,2,3$, if there is an uncountable set $S\subset X$ such that any pair of
distinct points from $S$ is distributionally scrambled of type $i$.\\

\section{Distributional chaos and factors}

We will show the existence of a system without any distributionally scrambled pair which is semiconjugated to a distributionally chaotic factor. This system is three-dimensional union of countably many homocentric cylinders with unit height and converging radius. First we state the following technical lemma about rotation on circle. Let $u\in\mathbb{S}$ and $v\in\mathbb{S}$ be determined by normed angles $\phi_u\in I$ and $\phi_v\in I$. These points rotate along the circle by different angles $r_u\in I$, respectively $r_v\in I$, i.e. \begin{equation}
\begin{split}
\phi_u&\mapsto(\phi_u+r_u) \text{ mod 1}\\ \phi_v &\mapsto(\phi_v+r_v)\text{ mod 1.}
\end{split}
\end{equation}
We denote the relative angle of rotation by $\Delta r=|r_u-r_v|$ and assume that the metric on $\mathbb{S}$ is $\rho(\alpha,\beta)=\min\{|\alpha-\beta|,1-|\alpha-\beta|\}$.
\paragraph{\bf Lemma 1}
\emph{For every number $\delta>0$ and every integer $p>\frac{2}{\Delta r}$, the following estimation holds:}
$$\frac{1}{p}\#\{0<i<p;\rho((\phi_u+ir_u)\text{ mod 1},(\phi_v+ir_v)\text{ mod 1})<\delta\}<3\delta.$$
\begin{proof}
Because $\rho((\phi_u+ir_u)\text{ mod 1},(\phi_v+ir_v)\text{ mod 1})=\rho(\phi_u,(\phi_v+i\Delta r)\text{ mod 1})$, it is sufficinet to show
$$\frac{1}{p}\#\{0<i<p;\rho(\phi_u,(\phi_v+i\Delta r)\text{ mod 1})<\delta\}<3\delta.$$
The expression $[p\cdot \Delta r]$ determines the number of turns that
 the point $v$ makes along the circle by rotation through the angle $\Delta r$ after $p$ iterations,
where $[x]$ denotes the integer part of $x$. The maximal number of iterations $i$ during one turn, for which $\rho(\phi_u,(\phi_v+i\Delta r)\text{ mod 1})<\delta$, is $\frac{2\delta}{\Delta r}.$ It follows 
$$\frac{1}{p}\#\{0<i<p;\rho(\phi_u,(\phi_v+i\Delta r)\text{ mod 1})<\delta\}<\frac{1}{p}([p\cdot \Delta r]\frac{2\delta}{\Delta r}+\frac{2\delta}{\Delta r})<2\delta+\frac{2\delta}{p\Delta r}.$$
Because $p>\frac{2}{\Delta r}$, we can estimate the second term by $\delta$, i.e.  $\frac{2\delta}{p\Delta r}<\delta$.

\end{proof}
\paragraph{\bf Theorem 1}
\emph{There exists a DC1 dynamical system $(Y, f)$ which is semiconjugated to an extension $(X,F)$ which possess no distributionally scrambled pair (of type 1 or 2).}

\begin{proof}
The space $X$ is defined $$X=\Big(\big \{[(2-\frac{1}{k})\cos2\pi\phi,(2-\frac{1}{k})\sin2\pi\phi]: k\in\mathbb{N},\phi\in I\big \}\bigcup\big \{[2\cos2\pi\phi,2\sin2\pi\phi]:\phi\in I\big \}\Big )\times I,$$ where $I$ is the unit interval. Each point $u=[r_u\cos2\pi\phi_u,r_u\sin2\pi\phi_u,z_u]$ in X is determined by its angle $\phi_u\in I$, radius $r_u\in\{2-\frac{1}{k}:k\in\mathbb{N}\}\cup\{2\}$ and height $z_u\in I$.\\
The space is endowed with max-metric \begin{equation}\label{metric}d(u,v)=\max\{|r_u-r_v|,|z_u-z_v|,\rho(\phi_u,\phi_v)\},\end{equation} \\where $\rho(\phi_u,\phi_v)=\min\{|\phi_u-\phi_v|,1-|\phi_u-\phi_v|\}$. 
We define the mapping $F:X\rightarrow X$ as identity on the limit cylinder,
$$[2\cos2\pi\phi,2\sin2\pi\phi,z]\mapsto[2\cos2\pi\phi,2\sin2\pi\phi,z],$$
and as a composition of rotation and continuous mapping $g$ on inner cylinders,
$$[(2-\frac{1}{k})\cos2\pi\phi,(2-\frac{1}{k})\sin2\pi\phi,z]\mapsto[(2-\frac{1}{k+1})\cos2\pi(\phi+\Psi(k,z)),(2-\frac{1}{k+1})\sin2\pi(\phi+\Psi(k,z)),g_{k}(z)].$$
To  define $g_k:I\rightarrow I$ and $\Psi:\mathbb{N}\times I\rightarrow I$, let $\{r_i\}_{i=1}^{\infty}=\mathbb{Q}|_{(0,1)}$ be a sequence of all rationals in $(0,1)$, and  $m_1<m_2<m_3<\ldots$ an increasing sequence of integers which we specify later. Then
\begin{equation}\label{definition_g}g_k = \left\{
  \begin{array}{l l}
    h_l & \quad \text{if $m_{3l+1}\leq k< m_{3l+2}$}\\
    Id&\quad \text {if $m_{3l+2}\leq k< m_{3l+3}$}\qquad k,l\in\mathbb{N}_0\\
    h_l^{-1}& \quad \text{if $m_{3l+3}\leq k< m_{3l+4}$}\\
    
  \end{array} \right.
  \end{equation}
  where $h_l:I\rightarrow I$ is a continuous strictly increasing mapping with three fixed points $0,1,r_l$ and $$\lim_{l\to\infty}||h_l-Id||=0; \quad h_l(x)<x\text{ for }x\in(0,r_l); \quad h_l(x)>x\text{ for }x\in(r_l,1).$$
  The sequence $\{m_i\}_{i=1}^{\infty}$ is defined in the following way:$$m_{3l+2}-m_{3l+1}=m_{3l+4}-m_{3l+3}=n_l,$$ where $n_l$ is integer satisfying
  \begin{equation}\label{inclusion}h_l^{n_l}([0,r_l-\frac{1}{l}])\subset[0,\frac{1}{l}) \quad\wedge \quad h_l^{n_l}([r_l+\frac{1}{l},1])\subset(1-\frac{1}{l},1],\end{equation} and simultaneously $\{m_i\}_{i=1}^{\infty}$ can be chosen such that \begin{equation}\label{min}m_{3l+3}-m_{3l+2}>\frac{2l}{\epsilon_l},\qquad\text{ where }\epsilon_l=\min\{h_l^{n_l}(\frac{1}{l}),1-h_l^{n_l}(1-\frac{1}{l})\},\end{equation} 
  
  \begin{equation}\label{distr}
  \begin{split}
 \lim_{l\to\infty}\frac{m_{3l+1}}{m_{3l+2}}&= \lim_{l\to\infty}\frac{m_{3l+3}}{m_{3l+4}}=1,\\
  \lim_{l\to\infty}\frac{m_{3l+2}}{m_{3l+3}}&=0.
    \end{split}\end{equation}
The angle of rotation $\Psi:\mathbb{N}\times I\rightarrow I$ is defined as
$$\Psi(k,z) = \left\{
  \begin{array}{l l}
    z&\quad \text{if $1\leq k<m_4$}\\
    z/l& \quad \text{if $m_{3l+1}\leq k< m_{3l+4}$}\qquad l\in\mathbb{N}.

  \end{array} \right.$$
  The factor space $Y$ is simply $X$ with fixed $\phi=0$, i.e. for each point $y\in Y,$
  $$y=[2-\frac{1}{k},0,z]\quad\text{or}\quad[2,0,z],\quad k\in \mathbb{N},z\in I.$$
  To simplify the notation, we skip the second zero coordinate and treat $Y$ as a two-dimensional space. The space $Y$ is union of converging sequence of unit fibers and the limit fiber,
  $$Y=\big\{2-\frac{1}{k}:k\in\mathbb{N}\big\}\times I\cup \{2\}\times I.$$
  Then the system $(X,F)$ is semiconjugated with skew-product map $f:Y\rightarrow Y$, which is identity on the limit fiber,
  $$[2,z]\mapsto[2,z],$$
  and which is $g_k$ on inner fibers,
  $$[2-\frac{1}{k},z]\mapsto[2-\frac{1}{k+1},g_k(z)],\qquad k\in\mathbb{N}.$$

\emph{I. The factor system (Y,f) is DC1.}\\
We show that set $S=\{1\}\times I$ is a distributionally scrambled set, i.e. for any pair of distinct points $(u,v)\in S^2$,\begin{equation}\label{factr}\Phi^*_{(u, v)}\equiv 1 \mbox{  and  } \Phi_{(u, v)}(\epsilon)=0, \mbox{  where  }\epsilon<1.
\end{equation}
Since $\{r_i\}_{i=1}^{\infty}$ is dense in $I$ and by (\ref{inclusion}), we can find a sequence $\{s_k\}_{k=1}^{\infty}$ such that $d(f^i(u),f^i(v))<\frac{1}{s_k}$, for $m_{3s_k+2}\leq i<m_{3s_k+3}$, and therefore, by (\ref{distr}), $\Phi^*_{(u, v)}\equiv 1$. Suppose $u^2>v^2$, where $x^2$ denotes the second coordinate of a point $x$. We can find another subsequence $\{q_k\}_{k=1}^{\infty}$ such that $d(f^i(u), f^i([1,1]))<\frac{1}{q_k}$ and simultaneously $d(f^i(v), f^i([1,0]))<\frac{1}{q_k}$, for $m_{3q_k+2}\leq i<m_{3q_k+3}$. Since $f$ preserves the distance between the endpoints of any fiber, $d(f^i([1,1]),f^i([1,0]))=1$, for $i\geq 0$, we can conclude, by (\ref{distr}), $\Phi_{(u, v)}(\epsilon)=0, \mbox{  for any  }\epsilon<1.$\\

\emph{II. $(X,F)$ has no distributionally scrambled pair}\\
We claim $\Phi^*_{(u,v)}<1$ for any pair of distinct points in $X$. Let $X_0$ be the limit cylinder $X_0=\{[2\cos2\pi\phi,2\sin2\pi\phi]:\phi\in I\}\times I$ and $\widetilde X=X\setminus X_0.$ Consider 4 possible cases:\\

a) $(u,v)\in\widetilde X$ with $z_u=z_v=z, \,k_u=k_v=k,\,\phi_u\neq\phi_v.$\\
The angle of rotation is the same for both $u$ and $v$, $\Psi(k_u,z_u)=\Psi(k_v,z_v)=\Psi(k,z)$, hence, by (\ref{metric}),
$$d(F(u),F(v))=\rho(\phi_u+\Psi(k,z),\phi_v+\Psi(k,z))=\rho(\phi_u,\phi_v)=d(u,v).$$ $F$ is isometric in this case and $\Phi^*_{(u,v)}\neq1$.\\

b)  $(u,v)\in\widetilde X$ with $z_u\neq z_v, \,k_u=k_v=k,\,\phi_u\neq\phi_v.$\\
Without loss of generality suppose $k=1$ (otherwise consider the pre images $(F^{-k}(u),F^{-k}(v))$) and let $L$ be an integer such that $|z_u-z_v|>\frac{1}{L}$. It is sufficient to show that there is $0<\delta<\frac{1}{3}$, for which $$\frac{1}{m_{3l+3}-m_{3l+2}}\#\{m_{3l+2}<i<m_{3l+3};d(F^i(u),F^i(v))<\delta\}<3\delta, \quad \text{for any} \quad L\leq l.$$
Since $d$ is max-metric, it is sufficient to prove $$\frac{1}{m_{3l+3}-m_{3l+2}}\#\{m_{3l+2}<i<m_{3l+3};\rho(\phi_{F^i(u)},\phi_{F^i(v)})<\delta\}<3\delta.$$
Since $|h_L^{n_L}(z_u)-h_L^{n_L}(z_v)|>\epsilon_L$ (see (\ref{definition_g}) and definition of $\epsilon_L$ in (\ref{min})), and $|h_L^{n_L}(z_u)-h_L^{n_L}(z_u)|$ is the minimal distance between trajectories of $u$ and $v$ between times $m_{3L+1}$ and $m_{3L+4}$, it follows 
\begin{equation}\label{min}\min_{3L+1<k\leq 3L+4}|g_k\circ g_{k-1}\circ\ldots\circ g_{3L+1}(z_u)-g_{k}\circ g_{k-1}\circ\ldots\circ g_{3L+1}(z_v)|>\epsilon_L.\end{equation}
Denote the relative angle of rotation of points with height $z_u$ and $z_v$ in the $k$-the cylinder by $\Delta \Psi_k(z_u,z_v)=|\Psi(k,z_u)-\Psi(k,z_v)|=\frac{|z_u-z_v|}{L}$, 
for $m_{3L+1}\leq k<m_{3L+4}.$ By (\ref{min}),
$$\Delta \Psi_k(g_{k}\circ g_{k-1}\circ\ldots\circ g_{3L+1}(z_u),g_{k}\circ g_{k-1}\circ\ldots\circ g_{3L+1}(z_v))>\frac{\epsilon_L}{L},\text{ for } m_{3L+1}\leq k<m_{3L+4}.$$
Since $m_{3L+3}-m_{3L+2}>\frac{2L}{\epsilon_L}$, we can use Lemma 1 and conclude, for any $\delta>0$, $$\frac{1}{m_{3L+3}-m_{3L+2}}\#\{m_{3L+2}<i<m_{3L+3};\rho(\phi_{F^i(u)},\phi_{F^i(v)})<\delta\}<3\delta.$$
We obtain the result for any $l>L$ using the same argument, since for every $l>L$, $|z_u-z_v|>\frac{1}{l}$.\\

c) $(u,v)\in\widetilde X$ with $z_u\neq z_v, \,k_u\neq k_v,\,\phi_u\neq\phi_v.$\\
Without loss of generality suppose $k_u=1$ and $k_v=p$. If $|z_u-z_v|>\frac{1}{L}$, then by case b) $$\#\{m_{4L+2}+p<i<m_{4L+3}-p;\rho(\phi_{F^i(u)},\phi_{F^i(v)})<\delta\}<3\delta\cdot (m_{3L+3}-m_{3L+2})$$
and hence
$$\frac{1}{m_{3L+3}-m_{3L+2}}\#\{m_{3L+2}<i<m_{3L+3};\rho(\phi_{F^i(u)},\phi_{F^i(v)})<\delta\}<3\delta+\frac{2p}{m_{3L+3}-m_{3L+2}}<1,$$ for sufficiently large $L$.

d) $u\in\widetilde X$ and $v\in X_0$\\
Since $v\in X_0$ is fixed and $\phi_v=\phi_{F(v)}$, we can find another point in $\widetilde X$, $w=[(2-\frac{1}{k_u})\cos2\pi\phi_v,(2-\frac{1}{k_u})\sin2\pi\phi_v,0]$, which is also fixed under rotation. Therefore $$\rho(\phi_{F(u)},\phi_{F(v)})=\rho(\phi_{F(u)},\phi_{F(w)})$$ and we can apply case b) or c) to investigate the pair $(u,w)$ instead of $(u,v)$.
\end{proof}
\paragraph{\bf Remark}
Notice that the upper distributional function for the extension remains positive, $\Phi^*_{(u,v)}>0$, for any pair of distinct points in $S\times \{1\}$. By (\ref{factr}), $\Phi_{(u,v)}<\Phi^*_{(u,v)},$ hence the system $(X,F)$ is distributionally chaotic of type 3. This fact implies an open question:
\emph{Is there a DC3 system wich is semiconjugated to an extension without any distributionally scrambled pairs of type 3?}

\paragraph{\bf Acknowledgment}
I sincerely thank my supervisor, Professor Jaroslav Sm\' ital, for valuable guidance. I~am grateful for his constant support and help.\\

\end{document}